\documentclass[12pt]{article}
\usepackage[latin1]{inputenc}

\usepackage{amsmath,amssymb,amstext,amsthm,amsfonts}
\usepackage{url}
\usepackage{color,graphicx,graphics,xcolor}
\usepackage{hyperref}
\usepackage{mathrsfs}
\usepackage{listings}

\usepackage{upref,amsxtra,exscale}
\usepackage{cite}
\usepackage{hyperref}
\usepackage {comment}

\usepackage{chngcntr}
\counterwithin*{equation}{section}

\newtheorem{thm}{Theorem}[section]
\newtheorem{defi}[thm]{Definition}

\newtheorem{pro}[thm]{Proposition}
\newtheorem{cor}[thm]{Corollary}
\newtheorem{lem}[thm]{Lemma}
\newtheorem{rem}[thm]{Remark}

\newtheorem{example}[thm]{Example}

\newcommand{\R}{\mathbb{R}}
\newcommand{\N}{\mathbb{N}}

\def\qq#1{\qquad \mbox{#1}\quad}
\newcommand{\D}{\displaystyle}

\newcommand{\al}{\alpha}
\newcommand{\be}{\beta}
\newcommand{\De}{\Delta}
\newcommand{\de}{\delta}

\newcommand{\e}{\varepsilon}

\newcommand{\la}{\lambda}

\newcommand{\om}{\omega}
\newcommand{\Om}{\Omega}
\newcommand{\Omb}{\overline{\Om}}
\newcommand{\p}{\partial}

\newcommand{\te}{\theta}

\newcommand{\vf}{\varphi}

\title{On the smoothness of weak solutions to subcritical semilinear elliptic equations in any dimension}
\author{Rosa Pardo
\footnote{Rosa Pardo, Universidad Complutense de Madrid, 28.040 Madrid, Spain, \url{rpardo@ucm.es}}
}
\date{}

\begin{document}
\maketitle
\begin{abstract}
Let us consider a semilinear boundary value problem
$ - \Delta u= f(x,u),$ in $\Omega,$ with Dirichlet boundary conditions, where $\Omega \subset \R ^N $,  $N> 2,$ is a bounded smooth domain.	
We provide sufficient conditions guarantying that   semi-stable weak positive solutions to subcritical semilinear elliptic equations are smooth in any dimension, and as a consequence, classical solutions. By a subcritical nonlinearity we mean $f(x,s)/s^\frac{N+2}{N-2} \to 0$ as $s\to\infty$, including non-power nonlinearities, and enlarging the class of subcritical nonlinearities, which is usually reserved for power like nonlinearities. 

\medskip

\medskip


\end{abstract}

\section{Introduction}
Let us consider the following semilinear boundary value problem
\begin{equation}
\label{eq:ell:pb} 
- \Delta u= f(x,u), \quad\mbox{in }\Om, \qquad  u= 0,\quad \mbox{on }\p \Om,
\end{equation}
where $\Om \subset \R^N $,  $N> 2,$ is a bounded, connected open subset, with $C^{2,\al}$ boundary $\p \Om$, 
and the non-linearity $f:\bar{\Om}\times\R\to \R$ is continuous and subcritical. Let $2^*=\frac{2N}{N-2}\ $ be the critical Sobolev exponent, by a {\it subcritical} nonlinearity we mean \begin{equation}\label{f:sub}
\lim_{s \to +\infty} \ \dfrac{\displaystyle\max_{x\in\Omb }f(x,s)}{s^{2^* -1 }}=0.
\end{equation} 
Usually the term subcritical nonlinearity is reserved for power like nonlinearities. Our analysis shows that nonlinearities satisfying \eqref{f:sub},  widen the  class of subcritical nonlinearities  sharing with power like nonlinearities properties such as $L^\infty$ a priori estimates, 
(see Theorem \ref{th:apriori:cnys} and Theorem \ref{th:apriori:cnys:k}), 
or regularity of semi-stable weak positive solutions.
 (see Theorem \ref{th:incr:convex}), and Theorem \ref{th:conv:seq}).
Our definition of a subcritical non-linearity includes  nonlinearities such as
$$
f^{(1)}(s):=\frac{(1+s)^{2^*-1}}{\big[\log(e+s)\big]^\be},\qq{or} f^{(2)}(s):= \frac{(1+s)^{2^*-1}}{\Big[\log\big[e+\log(1+s)\big]\Big]^\be},
$$
for any $\be>0$. 

\medskip

We focus in contributing to the problem of regularity of {\it weak solution} in the class of subcritical generalized problems, for any dimension $N>2$.

\smallskip

By a {\it solution} we mean a weak solution $u \in H_0^{1}(\Om)$  such that $f(x,u) \in L_{loc}^1(\Om)$, and  
\begin{equation}\label{weak:sol}
	\int_{\Om}\nabla u\nabla \vf=\int_{\Om} f(x,u)\vf, \quad\quad \forall \vf\in C_c^\infty(\Om).
\end{equation}
Let $u$ be a solution to \eqref{eq:ell:pb}. We will say that $u$ is  {\it semi-stable} if $f_s(\cdot,u)\in L^1_{loc}(\Om)$ and 
\begin{equation}\label{eq:semi:stab}
	\int_\Om|\nabla\vf|^2 \, dx\ge \int_\Om f_s(x,u)\vf^2 \, dx  , \qq{for all} \vf \in C_c^\infty(\Om),
\end{equation}
where $f_s:=\frac{\p f}{\p s}.$

Cabr\'e, Figalli, Ros-Oton, and Serra analyze the regularity  of  semi-stable\footnote{They call it stable solutions, see\cite[Definition 1.1]{Cabre_Figalli_Ros-Oton_Serra} } solutions with a nonlinearity $f=f(s)$ positive, non-decreasing,  convex, and such that $f(s)/s\to\infty$ as $t\to\infty$, and they conclude that semi-stable weak  solutions  in $H^1(\Om)$ are smooth up to dimension $N\le 9$, for domains of class $C^3$, see  \cite[Corollary 1.6]{Cabre_Figalli_Ros-Oton_Serra}. In their arguments, it is crucial  to assume $f$ to be convex,  non-decreasing and non-negative.

\medskip

Our aim  is to show  that, in addition to dimension, subcriticality   is another barrier dividing smoothness of  semi-stable solutions in $H^1(\Om)$. We show sufficient conditions guarantying that any set of positive semi-stable   solutions is uniformly $L^\infty (\Om)$ a priori bounded, and so they are classical solutions.

\bigskip

In order to prove our result, we first estimate the $L^\infty$-norm of {\it weak  solutions} in $H_0^{1}(\Om) \cap L^{\infty} (\Om)$,\footnote{\label{foot:reg} According to elliptic regularity, if $f$ is continuous in both variables, then $u$ is a {\it strong solution} in $W^{2,p}(\Om) \cap W_0^{1,p}(\Om)$, and by Sobolev embeddings, $u\in C^{1,\be}(\Omb)$ for any $\be<1$.} in terms of the $L^\frac{2N}{N+2}$-norm of $f$, see Theorem  \ref{th:apriori:cnys}. As we shall see below, the $L^\infty $ a priori bound of solutions requires nor convexity neither monotonicity of the nonlinearity. Moreover, we can prove that solutions are universally bounded  in terms only of their $L^{2^*}$- norm, with a constant  independent of  the solution, and surprisingly, independent of $f$ for non-decreasing nonlinearities in a neighborhood of infinity,  see Corollary \ref{cor:apriori:cnys}. In addition, this result holds for positive, negative and changing sign solutions. 

Secondly, we will approach weak solutions by smooth ones, see for instance \cite[Theorem 3]{Brezis_Cazenave_Martel_Ramiandrisoa} and \cite[Theorem 3.2.1 and
Corollary 3.2.1]{Dupaigne}. With this in mind, we work on sequences of BVP. More specifically, given a sequence of nonlinearities $f_k$, and the corresponding sequence of BVP, we provide sufficient conditions guarantying that any set of  solutions in $H_0^{1}(\Om) \cap L^{\infty} (\Om)$ is uniformly $L^\infty (\Om)$ a priori bounded, see Theorem \ref{th:apriori:cnys:k}.  

As an application, we next state another main result,  concerning the global regularity of semi-stable positive solutions in any dimension, when the non-linearity is subcritical, convex and non-decreasing. Assuming that $u^*\ge 0$ is a semi-stable  weak solution to \eqref{eq:ell:pb}, we build a sequence of  non-negative solutions in $H_0^{1}(\Om) \cap L^{\infty} (\Om)$ upper bounded by $u^*$.
The key point here is the uniform $L^\infty$ a-priori bounds for that sequence of solutions.
Thanks to that and to the elliptic regularity, we obtain a  subsequence, convergent to $\tilde{u}$ in $C^{1,\be}(\Omb)$ for any $\be<1$. More regularity on $f$ guaranties more regularity on $\tilde{u}$. The limit $\tilde{u}$ is clearly a solution to \eqref{eq:ell:pb}. To conclude that $u^*$ is a classical solution, we need to prove that in fact, the limit is $u^*$, and at this point, we use the convexity of $f$. As a consequence, weak solutions are classical solutions, see  Theorem \ref{th:incr:convex}.

Finally, we elude hypothesis on convexity in  Theorem \ref{th:conv:seq}, using monotonicity methods, and still giving sufficient conditions so that  weak  solutions  in $H^1(\Om)$ to subcritical elliptic equations are smooth in any dimension. We emphasize that this result holds for weak  solutions, not necessarily semi-stable.

\bigskip

Smoothness of semi-stable weak solutions is a very classical topic in elliptic equations, posed 
by Joseph and Lundgren in \cite{Joseph_Lundgren}. They work on particular nonlinearities $f=f(s)$, with $f(s):=e^s$ or $f(s):=(1+s)^p$.  They consider the following BVP depending on a multiplicative parameter $\la\in\R$,
\begin{equation}\label{pde:la} 
- \Delta u=\la f(u), \quad\mbox{in }\Om, \qquad  u= 0,\quad \mbox{on }\p \Om,
\end{equation}
and look for classical radial positive solutions in  the unit ball $B_1$.
Furthermore,  they study singular solutions 
as limit of classical solutions.

\smallskip

When $N>2$, $\la=2(N-2)$,  and $f(s):=e^s$, they obtain the explicit weak solution
\begin{equation}\label{sing:sol:exp}
u^{*}_1(x):=\log \frac{1}{|x|^2},
\end{equation}
see \cite[p. 262]{Joseph_Lundgren}. It can be seen that $u^{*}_1\in H_0^1(B_1)$, and that $u^{*}_1$ is a singular weak solution to \eqref{pde:la} in the unit ball. 

On the other hand, the Hardy inequality states that
\begin{equation}\label{ineq:Hardy}
\int_\Om |\nabla\vf|^2 \, dx\ge \left(\frac{N-2}{2}\right)^2
\int_\Om \frac{\vf^2}{|x|^2}\, dx
, \qq{for all} \vf \in C_0^1(\Om)
\end{equation}
when $N\ge 3$, and then $u^{*}_1$ is a singular semi-stable solution when $N\ge 10$. Observe that $f(s):=e^s$ is not a subcritical non-linearity.

\smallskip

When $N>2$, $f(s):=(1+s)^p$ with $p>\frac{N}{N-2}$, and $\la=\frac{2}{p-1} \big(N-\frac{2p}{p-1}\big)$,  they also found the explicit {\it $L^1$-weak solution} 
\begin{equation}\label{sing:sol:pol}
u^{*}_2(x):= \left(\frac{1}{|x|}\right)^\frac{2}{p-1}-1,
\end{equation}
see \cite[(III.a)]{Joseph_Lundgren}. 
We will say that a function $u$ is a  $L^1$-{\it weak solution} to \eqref{eq:ell:pb} if 	
$$
u\in L^1(\Om),\qquad f(\cdot,u)\, \de_\Om \in   L^1(\Om)
$$
where $\de_\Om (x):=dist(x,\p\Om)$ is the distance function with respect to the boundary, and 
\begin{equation}\label{eq:L1:sol}
\int_\Om \Big(u  \Delta \varphi +  f(x,u)\varphi\Big) \, dx=0, \qq{for all} \varphi \in C^2(\Omb ),\quad \varphi\big|_{\p\Om}=0. 
\end{equation}
It holds that $u^{*}_2\in W_0^{1,\frac{N}{N-1}}(B_1)$ for $p>\frac{N}{N-2}$, and $u^{*}_2$ is a singular $L^1$-weak
solution to \eqref{pde:la} on the unit ball. Moreover $u^{*}_2\in H_0^1(B_1)$ only when $p>\frac{N+2}{N-2}$.

Since the Hardy inequality \eqref{ineq:Hardy}, it can be checked that $u^{*}_2$ is a semi-stable solution if $\big(\frac{N-2}{2}\big)^2\ge \frac{2p}{p-1} \left(N-\frac{2p}{p-1}\right)$. 
Note that $f(s):=(1+s)^p$ is a subcritical non-linearity whenever $p<\frac{N+2}{N-2}$. In the subcritical range, $u^{*}_2$ is a semi-stable solution for $p\le\frac{N+2\sqrt{N-1}}{N-4+2\sqrt{N-1}}$, and  $u^{*}_2\in W_0^{1,\frac{N}{N-1}}(B_1)$, so  $u^{*}_2$ is a singular $L^1$-weak
solution, not in $H^1$.

\medskip

Those examples for radially symmetric solutions to BVP's on spherical domains show that the existence of  singular solution(s) in $H_0^1(\Om)$ is  not only  related with the dimension, but also with the sub-critical, {\it critical}, or {\it supercritical}  nature of the non-linearity.
By a {\it critical} ({\it supercritical}) non-linearity we mean $f(x,s)=O\big(s^\frac{N+2}{N-2}\big)$,  $\big(f(x,s)\big/ s^\frac{N+2}{N-2}\to\infty\big)$ respectively, as $s\to\infty$.

\medskip

It is natural to ask for the extent of these results on singular positive solutions, over more general nonlinearities and  non-spherical domains.
The regularity of semi-stable solutions to semilinear elliptic equations, is initiated in \cite{Joseph_Lundgren} with the explicit examples already mentioned,  continued by Keener and Keller \cite{Keener_Keller}, and by Crandall and Rabinowitz  in \cite{C-R_1975}, and rising a huge literature on the topic, see the monograph \cite{Dupaigne} for an extensive list of results and references.
Crandall and Rabinowitz consider a nonlinearity $f\in C^3$, positive, non-decreasing, convex,  and superlinear at infinity. They state  that if $N<10$ and the following limit exists
$$
\lim_{s \to \infty}\,\frac{f(s)\, f''(s)}{\big(f_s(x,s)\big)^2}:= L\ge 0,
$$
then  $H_0^1(\Om)$ semi-stable solutions to \eqref{pde:la} are bounded.
Brezis and Vázquez study singular  {\it $L^1$-weak solutions}, unbounded in $L^\infty$, for nonlinearities $f\in C^2$, positive, non-decreasing, convex,  and superlinear at infinity, see \cite{Brezis_Vazquez}. 
When $f(s)=s^\frac{N}{N-2}$, Pacard in \cite{Pacard} prove the existence of positive $L^1$-weak solutions with prescribed singular set.
R\'{e}bai in \cite{Rebai_96} study the existence of positive $L^1$-weak solutions  which are singular either at exactly $N$ points, for $N\ge 2$, or on a prescribed $(N - m)$-dimensional compact submanifold $\Sigma\subset\Om$ without boundary, with $N> m>2$, when $f(s)=s^p$ for $p>\frac{m}{m-2}$ and close to that number. In both cases,  those $L^1$-weak solutions are not in $H^1(\Om)$. Results on supercritical problems and their singular sets can be read in \cite{Guo-Li} and references therein.

\bigskip

This paper is organized in the following way. In Section \ref{sec:main} we state our main results.
In Section \ref{sec:prelim} we include some preliminaries  and known results. Section \ref{sec:proof} contains the proofs of  Theorem \ref{th:apriori:cnys} and Theorem \ref{th:apriori:cnys:k}. Section \ref{sec:proof2} is devoted to prove Theorem \ref{th:incr:convex} and Theorem \ref{th:conv:seq}.

\medskip

\section{Main results: Estimates of the $L^\infty$-norm  of the solutions and Regularity of semi-stable weak solutions}
\label{sec:main}
In this Section, we state our main results. We will do it in two parts. In the first part 
we estimate the $L^\infty$-norm of the solutions, in terms of the $L^\frac{2N}{N+2}$-norm of $f$, see Theorem \ref{th:apriori:cnys}. Also, 
given a sequence of nonlinearities $f_k$, and the corresponding sequence of BVP, we estimate the $L^\infty$-norm of solutions to the sequence of BVP in terms of the $L^\frac{2N}{N+2}$-norm of $f_k$, see  Theorem \ref{th:apriori:cnys:k}. This results hold for positive, negative and changing sign solutions. 

In the second part, we show that positive semi-stable weak solutions can be approximated by smooth ones,  when $f$ is convex, see Theorem \ref{th:incr:convex}. We also show that weak solutions can be approximated by smooth ones, when $f_s>\la_1$, see  Theorem \ref{th:conv:seq}.

\subsection{Part I. Estimates of the $L^\infty$-norm of the solutions}
We assume that the nonlinearity $f:\bar{\Om}\times\R\to\R$  is continuous in both variables and satisfy the following assumptions
\begin{enumerate}				
\item[{\rm (H1)}] $f$ is subcritical, that is  
$\ \D\lim_{|s| \to \infty} \ \dfrac{\displaystyle\max_{x\in\Omb }|f(x,s)|}{|s|^{2^* -1 }}
=0\ $ where $2^*=\frac{2N}{N-2}$ is the critical Sobolev exponent.

\item[\rm (H2)] there exists a uniform constant $c_0>0$ such that 
\begin{equation}\label{H2}
\limsup_{s\to +\infty}\ \dfrac{\displaystyle \max_{\Omb \times[-s,s]}\, |f|} {\displaystyle\max_{\Omb\times\{-s,s\}} |f|}\, \le \, c_0.
\footnote{Observe that in particular, if  $|f(x,\cdot)|$ is monotone for all $x\in\Om$, then (H2) is obviously satisfied with $c_0=1$.} 
\end{equation}

\item[\rm (H3)] there exists two constants $M_0>0$ and $s_0>0$ such that 
\begin{equation}\label{H3}
	\max_{x\in\Omb}|f(x,s)|>M_0 \qq{for} |s|>s_0.
\end{equation}	
\end{enumerate}

Let us define
\begin{equation}\label{def:h}
h(s):=\frac{|s|^{2^* -1 }}{\displaystyle 
\max_{\Omb\times\{-s,s\}} |f|}, \qq{for} |s|> s_0.
\end{equation}

\bigskip

Our first main results is the following  theorem. Let  $u\in H_0^{1}(\Om) \cap L^{\infty} (\Om)$ be any weak solution   to \eqref{eq:ell:pb}.
\footnote{Despite the regularity inherent to an $L^\infty$ bound, (see footnote {\footnotesize \ref{foot:reg}}), we keep the notation as above $\big($weak solution  in $H_0^{1}(\Om) \cap L^{\infty} (\Om)\big)$, in order to clarify which hypothesis are specifically involved in each  statement.}
Under hypothesis  (H1)-(H3), we establish an estimate for the function $h$ applied to the $L^\infty (\Om)$-norms of the solution, in terms of  the $L^\frac{2N}{N+2} (\Om)$ norm of $f$. 

From now on, $C$ denotes several constants that may change from line to line, and are independent of $u$. 

\begin{thm}
\label{th:apriori:cnys}  
Assume that $f:\Omb \times\R\to \R$ is a  
continuous function in both variables  satisfying  {\rm (H1)}-{\rm (H3)}. 

Then,   there exists a constant $C>0$ such that for any  $u\in H_0^{1}(\Om) \cap L^{\infty} (\Om)$ weak solution   to \eqref{eq:ell:pb}, the following holds:
\begin{enumerate}				
\item[{\rm (i)}] either  $\|u\|_{\infty}\le C$, where  $C$  is independent of the  solution $u$,
\item[{\rm (ii)}] or for any $\te\in(0,1]$
\begin{equation}\label{L:inf:h:0}
h\big(\|u\|_{\infty}\big)\le C\ \big(\|u\|_{2^*}\big)^{\frac{N+2}{N-2}\frac{2(1-\te)}{N-2\te}} \ \Big(\|f(\cdot,u) \|_{\frac{2N}{N+2}}\Big)^\frac{2(1+\te)}{N-2\te},
\end{equation}
where $h$ is defined by \eqref{def:h}, and $C$ depends only on $\Om,$ $\te$ and $N$ and  it is independent of the  solution $u$.

In particular (for  $\te=1$)
\begin{equation}\label{L:inf:h:3}
h\big(\|u\|_{\infty}\big)
\le C\,\Big(\|f(\cdot,u) \|_{\frac{2N}{N+2}} \Big)^\frac4{N-2},
\end{equation}
where  $C$ depends only on $\Om,$ and $N$ and it is independent of  $u$.
\end{enumerate}
\end{thm}

As as immediate corollary, we prove that  any sequence of solutions uniformly bounded in the $ L^{2^*} (\Om)$-norm,  is also uniformly bounded in the $ L^{\infty} (\Om)$-norm.

\begin{cor}\label{cor:apriori:cnys}
Let $f:\Omb \times\R\to \R$ be a  continuous function in both variables satisfying {\rm (H1)}--{\rm (H3)}. 

Let  $\{u_k\}\subset H_0^{1}(\Om) \cap L^{\infty} (\Om)$ be any sequence of solutions to \eqref{eq:ell:pb} such that there exists a constant $C_0>0$ satisfying
\begin{equation}\label{L:2*:k:0}
\|u_k\|_{2^*} \le C_0.
\end{equation}

Then, there exists a constant $C>0$ such that
\begin{equation}\label{L:inf:k:0}
\|u_k\|_{\infty} \le C.
\end{equation} 
\end{cor}

\begin{proof}
We reason by contradiction, assuming that \eqref{L:inf:k:0} does not hold.
Indeed, by  subcriticallity, and \eqref{L:2*:k:0}, 
\begin{align*}
\Big(\|f(\cdot,u_k) \|_{\frac{2N}{N+2}}\Big)^{\frac{2N}{N+2}}
&\le C\left(1+\int_\Om |u_k|^{2^*}\, dx\right)\le C.
\end{align*}

Now part (ii) of the Theorem \ref{th:apriori:cnys} implies that 
\begin{equation}\label{h:bdd:uk:0}
h\big(\|u_k\|_{\infty}\big)
\le C.
\end{equation}
From \eqref{def:h}  and hypothesis (H1), for any $\e_0>0$ there exists $s_0>0$ such that  $h(s)\ge 1/\e_0$ for any $s\ge s_0$. This, joint with \eqref{h:bdd:uk:0} ends the proof.	
\end{proof}

\subsubsection{Estimates of the $L^\infty$-norm of  solutions to sequences of BVP}

Next, we state our second main result. It  concerns sequences of subcritical BVP.

Let us now consider a sequence $f_k:\Omb \times\R\to \R$ of continuous functions in both variables satisfying the following conditions:
\begin{enumerate}
	\item[\rm (H2)$_k$]  There exists a uniform constant $c_1>0$ such that 
	\begin{equation}\label{H2:k}
		\limsup_{s\to +\infty}\ \dfrac{\displaystyle \max_{\Omb \times[-s,s]}\, |f_k|} {\displaystyle\max_{x\in\Omb}|f_k(x,s)|}\, \le \, c_1.
	\end{equation}
	\item[\rm (H3)$_k$] there exists two constants $M_0>0$ and $s_0>0$ such that 
	\begin{equation}\label{H3:k}
		\max_{x\in\Omb}|f_k(x,s)|>M_0 \qq{for} |s|>s_0.
	\end{equation}	
\end{enumerate}

Let us also consider the corresponding sequence of  elliptic equations:
\begin{equation}\label{eq:ell:pb:k}
\left\{ 
\begin{array}{rcll}
- \De u &=& f_k(x,u)  &\qquad\text{in}\ \Om,\\
u &=& 0 &\qquad\text{on}\ \p \Om.
\end{array}
\right.
\end{equation}
For each $k\in\N$, let $u_k\in H_0^{1}(\Om) \cap L^{\infty} (\Om)$ be a  solution to \eqref{eq:ell:pb:k}$_k$. Consider a sequence $\{u_k\} $ of those solutions. 

\medskip

In the following Theorem, we state sufficient conditions for having a uniform $L^{\infty}$  estimate for sequences of solutions $\{u_k\}\subset H_0^{1}(\Om) \cap L^{\infty} (\Om)$ to  \eqref{eq:ell:pb:k}$_k$.

\begin{thm}
\label{th:apriori:cnys:k}  
Assume that $f_k:\Omb \times\R\to \R$ is a  sequence of
continuous functions in both variables  satisfying   {\rm (H1)} for $f=f_k$, and  {\rm (H2)$_k$}-{\rm (H3)$_k$}. 

Then,   there exists a constant $C>0$ such that for any sequence $\{u_k\}\subset H_0^{1}(\Om) \cap L^{\infty} (\Om)$ of solutions   to \eqref{eq:ell:pb:k}$_k$, the following holds:
\begin{enumerate}				
\item[{\rm (i)}] either $\|u_k\|_{\infty}\le C$,
\item[{\rm (ii)}] or for any $\te\in(0,1]$
\begin{equation}\label{L:inf:hk:0}
h_k\big(\|u_k\|_{\infty}\big)\le C\ \big(\|u_k\|_{2^*}\big)^{\frac{N+2}{N-2}\frac{2(1-\te)}{N-2\te}} \ \Big(\|f_k(\cdot,u_k) \|_{\frac{2N}{N+2}}\Big)^\frac{2(1+\te)}{N-2\te},
\end{equation}
where $h_k$ is defined by \eqref{def:h} for $f=f_k$, and $C=C(\Om,\te,N)$ and it is independent of $k$. 

In particular (for  $\te=1$)
\begin{equation}\label{L:inf:hk:3}
h_k\big(\|u_k\|_{\infty}\big)
\le C\,\Big(\|f_k(\cdot,u_k) \|_{\frac{2N}{N+2}} \Big)^\frac4{N-2},
\end{equation}
where  $C=C(\Om,N)$ and it is independent of $k$.
\end{enumerate}
\end{thm}

\medskip

The arguments are inspired in the equivalence between uniform $L^{2^*}(\Om)$ {\it a priori} bound and uniform $L^\infty (\Om)$ {\it a priori} bound for solutions to subcritical elliptic  equations, see \cite[Theorem 1.2]{Castro_Mavinga_Pardo} for the semilinear case and $f=f(u)$, and \cite[Theorem 1.3]{Mavinga_Pardo_MJM} for the $p$-laplacian  and $f=f(x,u)$. Related results can be found in 
\cite{Castro_Pardo_RMC, Damascelli_Pardo, Castro_Pardo_DCDS, Mavinga_Pardo17_JMAA, Castro_Pardo_DCDS-A, Pardo19a, Pardo19b}. 

\bigskip

\subsection{Part II. Approximation of some weak solutions by sequences of classical solutions}

In this second Part, we apply the above results on $L^\infty$ a priori bounds, to positive semi-stable weak solutions. We consider non-negative functions. The Maximum Principle ensures that solutions are now non-negative.

We assume that the nonlinearity  $\, f:\Omb \times [0,\infty) \to [0,\infty)$ satisfy some of the following assumptions: 
\begin{enumerate}				
\item[\rm (H4)] There exist a constant  $c_0>1$, such that 
\begin{equation}\label{H4}
\D\liminf_{s \to +\infty}\ \frac{sf_s(x,s)}{f(x,s)}
\ge c_0>1 ,\qq{where} f_s(x,s):=\frac{\p f}{\p s}(x,s).
\end{equation}

\item[\rm (H5)] There exist a positive constant  $c_1$, such that 
$$ 
\inf_{\Omb\times\R}\ f_s(x,s)\ge c_1>\la_1,\quad \forall (x,s)\in\Omb\times\R ,
$$
where $\lambda_1:=\la_1(- \De;\Om )$ is the first eigenvalue of $- \Delta$ acting on $H^1_0(\Om)$
\end{enumerate}

\bigskip

Our first main result is the following Theorem, showing that  any positive  weak solution $u^{*}\in H_0^1(\Om)$ of a subcritical elliptic problem, with $f$ non-negative, non-decreasing, convex and satisfying (H4), is in fact a classical solution.  
Roughly speaking, this result is known for  nonlinearities $f=f(u)$ (non-negative, non-decreasing, and convex, but not necessarily subcritical),   
when $N\le9$, see  
\cite[Corollary 1.6]{Cabre_Figalli_Ros-Oton_Serra}.

\begin{thm}\label{th:incr:convex}
Let  $\, f:\Omb \times[0,\infty)\to [0,\infty)$ be  continuous in both variables, and continuously derivable with respect to the second variable. Assume also that $\forall x\in\Om, \ f(x,\cdot)$ is  non-decreasing, and convex, that $f(\cdot,s)\in C^\al(\Omb)$ for all $s\ge 0$, and that $f$ satisfies {\rm (H1)}-{\rm (H4)}. 

Let $u^{*}\in H_0^1(\Om)$ be a semi-stable weak solution  to \eqref{eq:ell:pb}. 

\smallskip

If 
\begin{equation}\label{0:no:sol}
f(x,0)\gneq 0,
\end{equation}
then, 
$u^{*}\in C^{2,\al}(\Omb )$ is a classical solution.
\end{thm}

\medskip

A sketch of the proof is the following: a weak solution $u^{*}\in H_0^1(\Om)$ is the limit  of a curve of smooth solutions to  approximate equations. A uniform $L^{2^*}$ bound is the key.
The equivalence between a uniform $L^{2^*}(\Om)$ {\it a priori} bound and a uniform $L^\infty (\Om)$ {\it a priori} bound for weak solutions to subcritical elliptic  equations implies a uniform $L^\infty (\Om)$ {\it a priori} bound (see \cite[Theorem 1.2]{Castro_Mavinga_Pardo} for the semilinear case and $f=f(u)$, and \cite[Theorem 1.3]{Mavinga_Pardo_MJM} for the $p$-laplacian  and $f=f(x,u)$). 
By Schauder elliptic regularity, the sequence is uniformly bounded in $ C^{2,\al}(\Omb )$. By compactness and monotonicity, the approximate solutions converges to $\tilde{u}\le u^*$ in $C^{2,\be}(\Omb )$ for any $\be<\al$. If $u^{*}$ is semi-stable, 
hypothesis (H4) and  convexity   play  an overriding role to prove that in fact $\tilde{u}= u^*$.

\medskip

Our second main result focuses on proving regularity, excluding convexity.  We  use instead hypothesis (H5) and sub and supersolution methods, proving the existence of a minimal and a maximal solution and  that the sequence of maximal solutions converge to $u^*$. 
It is applicable to any positive weak solution, independently of its stability.

\begin{thm}\label{th:conv:seq}	
Let  $\, f:\Omb \times[0,\infty)\to [0,\infty)$ be  continuous in both variables, and continuously derivable with respect to the second variable. Assume also that  $f(\cdot,s)\in C^\al(\Omb)$ for all $s\ge 0$, and $f$ satisfies {\rm (H1)}-{\rm (H3)}, and	{\rm (H5)}. 

Let $u^{*}\in H_0^1(\Om)$ be a non-negative weak solution to \eqref{eq:ell:pb}.

If \eqref{0:no:sol} is satisfied, then  $u^{*}\in C^{2,\al}(\Omb )$ is a classical solution.
\end{thm}

\section{Preliminaries  and known results}
\label{sec:prelim}
Consider the Hilbert space $H_0^{1}(\Om):=\{u\in L^{2}(\Om):\nabla u\in L^2(\Om,\R^N)\}$ with its usual inner product and norm,
$$
\langle u,v\rangle:=\int_\Om\nabla u\cdot\nabla v,\qquad \|u\|:=\left(\int_\Om |\nabla u|^2\right)^{1/2}.
$$
We will denote by $\|\cdot\|_p$  the standard $L^p$-norm.

Let $\la_1:=\la_1(- \De;\Om )$, and let $\phi_1>0$ denote the corresponding eigenfunction, normalized in the $L^\infty$-norm. 

\medspace

By elliptic regularity,   $L^1$-weak solutions to \eqref{eq:ell:pb}, bounded in $L^\infty$ are strong or classical solutions.

\begin{pro}\label{pro:classical}
Assume that $\p\Omega$ is $C^{2,\al}$. 
Assume that $f:\Om\times\R\to\R$  is a continuous function in both variables. Let $u$ be a  $L^1$-weak solution to \eqref{eq:ell:pb}. 

If $u\,  \in  L^\infty(\Om),$ then the following holds:
\begin{enumerate}
\item[\rm (i)] $u$ is a strong solution in $W^{2,p}(\Om)\cap W_0^{1,p}(\Om)$ for any $p>1$, and	$u\in  C^{1,\be}(\overline{\Om})$ for any $\be<1$ satisfies
\begin{equation}\label{u:C1:be}
\|u\|_{C^{1,\be}(\overline{\Om})}\le C \|f(\cdot, u)\|_{L^\infty(\Om)}. 
\end{equation}

\item[\rm (ii)]	Moreover, if $f:\Om\times\R\to\R$  is such that for all $R>0$ there exists $L=L(R)>0$ satisfying
$$
|f(x,s)-f(y,t)|\le L\big(|x-y|^\al+|s-t|^\al\big),
$$
for all $s,t\in[-R,R],\ x,y\in\Om,$ then $u\in C^{2,\al}(\overline{\Om})$ is a classical solution and
\begin{equation}\label{u:C2:al}
	\|u\|_{C^{2,\al}(\overline{\Om})}\le C \|f(\cdot, u)\|_{C^\al(\Omb)}. 
\end{equation}
\end{enumerate}

\end{pro}
\begin{proof}
(i) Since $u\in L^\infty(\Om),$	and $f$ is continuous  in both variables, $f(\cdot, u)\in L^\infty(\Om).$ By Agmon-Douglis-Nirenberg elliptic regularity, $u\in W^{2,p}(\Om),$ for any $p>1$, moreover $\|u\|_{W^{2,p}(\Om)}\le C\|f(\cdot,u)\|_{L^\infty(\Om)}$. By Sobolev embeddings $u\in C^{1,\be}(\overline{\Om}),$ for any $\be<1$
and estimate \eqref{u:C1:be} holds. 

(ii) If $\p\Omega$ is $C^{2,\al}$, by Schauder elliptic regularity, $u\in C^{2,\al}(\overline{\Om})$ (see \cite[Theorem 6.8]{G-T}) and since \cite[Theorem 6.6]{G-T} with homogeneous Dirichlet boundary conditions, estimate \eqref{u:C2:al} holds. 
\end{proof}

\begin{lem}\label{lem:w:H}
Assume that $f$ is subcritical.
Then, for any $u$ weak solution to \eqref{eq:ell:pb} 
the following hold:
\begin{equation}\label{f:L}
f(\cdot,u)\in L^{\frac{2N}{N+2}}(\Om),
\end{equation}
and
\begin{equation}\label{eq:E-L:H}
\int_{\Om}\, \Big(\nabla u\cdot \nabla \varphi  -  f(x,u)\varphi\Big)\,  dx=0, \qq{for any}\varphi\in H_0^1(\Om).
\end{equation} 
\end{lem}
\begin{proof}
Using \eqref{f:sub}, and  Sobolev embeddings, for any $u\in H_0^1(\Om)$ there exists a constant $C>0$ such that
$$
\int_{\Om} |f(x,u)|^\frac{2N}{N+2}\,  dx
\le 
C\left(1+\int_{\Om}|u|^{2^*}\,  dx \right)
\le C\left(1+\|u\|_{H_0^1(\Om)}^{2^*}\right)<+\infty,
$$
hence \eqref{f:L} holds.

In addition, by Holder inequality 
$$
\int_{\Om}\,  f(x,u)\varphi\,  dx<+\infty \qq{for any}\varphi\in H_0^1(\Om),
$$ 
and by density,  \eqref{eq:E-L:H} holds.
\end{proof}

\section{Proof of Theorems \ref{th:apriori:cnys} and \ref{th:apriori:cnys:k}}
\label{sec:proof}

In that Section, we  prove Theorem \ref{th:apriori:cnys:k} for sequences of BVP, and sequences of solutions $\{u_k\}$ to \eqref{eq:ell:pb:k}. The proof of Theorem \ref{th:apriori:cnys} is a particular case of Theorem \ref{th:apriori:cnys:k}, applied to one particular BVP, \eqref{eq:ell:pb}, and we omit it.

\begin{proof}[Proof of Theorem \ref{th:apriori:cnys:k}]
Let   $\{u_k\}\subset H_0^{1}(\Om) \cap L^{\infty}(\Om)$ be a sequence of  weak solution to \eqref{eq:ell:pb:k}$_k$.  
If $\|u_k \|_{\infty}\le C,$ then 
(i)  holds.

Now, we  argue on the  contrary, assuming that $\|u_k \|_{\infty} \to + \infty$ as $k \to \infty$

Let $x_k \in \Om$ be such that
$$|u_k(x_k)| =  \max_{\Omb } |u_k| .$$
Choose  $R_k$ such that  
\begin{equation*}\label{Sobolev:87:k}
	|u_k (x)|\ge
	\frac12\ \|u_k \|_{\infty}\qquad\mbox{for any}\quad  x\in \overline{B}(x_k,R_k),
\end{equation*}  
and there exists $y\in\p B(x_k,R_k)$ such that
\begin{equation}\label{Sobolev:88:2:k}
	|u_k (y)|=\frac12\ \|u_k \|_{\infty}.
\end{equation}

\bigskip

{\it Step 1. $W^{2,q}$ estimates for $q\in(N/2,N)$.}

\medskip

\noindent Let us denote by 
\begin{equation}\label{def:M:k:f}
M_k:=\max_{\Omb\times \{-\|u_k\|_{\infty},\|u_k\|_{\infty}\} }\big| f_k\big| \ge 
C\max_{\Omb \times [-\|u_k\|_{\infty},\|u_k\|_{\infty}]}|f_k|,
\end{equation}
by hypothesis (H2)$_k$, see \eqref{H2:k}.

For any $q>\frac{2N}{N+2},$
\begin{eqnarray}\label{f:q:k:2}
\displaystyle \int_{\Om} \left|f_k\big(x,u_k(x)\big)\right|^q\, dx  
&\le&  \int_{\Om} \left|f_k\big(x,u_k(x)\big)\right|^{\frac{2N}{N+2}} \, 
\left|f_k\big(x,u_k(x)\big)\right|^{q-\frac{2N}{N+2}}\, dx \nonumber\\
&\le & C \Big(\|f_k(\cdot,u_k) \|_{\frac{2N}{N+2}}\Big)^\frac{2N}{N+2}\
M_k^{\ q-\frac{2N}{N+2}}.
\end{eqnarray}

Let us take $q$ in the interval  $(N/2,N).$ Combining elliptic regularity with Sobolev embedding,   we have that
\begin{equation}\label{ell:reg:k}
\|u_k \|_{W^{1,q^*} (\Om)}
\le   C\
\Big(\|f_k(\cdot,u_k) \|_{\frac{2N}{N+2}}\Big)^\frac{2N}{q(N+2)}\
M_k^{\ 1-\frac{2N}{q(N+2)}},
\end{equation}
where $1/q^*=1/q-1/N$; since $q>N/2$, then $q^*>N.$

\bigskip

{\it Step 2. A lower bound for the radius $R_k$.}

\medskip

\noindent Using  Morrey's Theorem, we have that  
\begin{equation}\label{Morrey:2:k}
|u_k(x_1)-u_k(x_2)|\le C |x_1-x_2|^{2-N/q}\|\nabla  u_k\|_{q^*},\quad \forall x_1,x_2\in \Om,
\end{equation}
where the constant $C$ depends only on $\Om,$ $q$ and $N$.
Hence, for all $x\in \overline{B}(x_1,R)\subset\Om$
\begin{equation}\label{Morrey:1:k}
|u_k(x)-u_k(x_1)|\le  C\ R^{2-\frac{N}{q}}\|\nabla  u_k\|_{q^*},
\end{equation}
for any $k$.
In particular, it follows that  for any $x\in \overline{B}(x_k,R_k),$
\begin{equation}\label{Sobolev:2:k}
|u_k (x)-u_k (x_k)|\le C\ (R_k)^{2-\frac{N}{q}}\
\Big(\|f_k(\cdot,u_k) \|_{\frac{2N}{N+2}}\Big)^\frac{2N}{q(N+2)}\
M_k^{\ 1-\frac{2N}{q(N+2)}}.
\end{equation}

Taking $x=y$ in the above inequality and from \eqref{Sobolev:88:2:k} we obtain
\begin{equation}\label{Sobolev:3:k}
C\ (R_k)^{2-\frac{N}{q}}\
\Big(\|f_k(\cdot,u_k) \|_{\frac{2N}{N+2}}\Big)^\frac{2N}{q(N+2)}\
M_k^{\ 1-\frac{2N}{q(N+2)}}
\ge \frac12 \|u_k \|_{\infty},
\end{equation}
which implies
\begin{equation}\label{Sobolev:4:k}
(R_k)^{2-\frac{N}{q}}\ge  \frac1{2C}\,
\frac{\|u_k \|_{\infty}\ }{\Big(\|f_k(\cdot,u_k) \|_{\frac{2N}{N+2}}\Big)^\frac{2N}{q(N+2)}\
M_k^{\ 1-\frac{2N}{q(N+2)}}}  ,
\end{equation}
or equivalently
\begin{equation}\label{Sobolev:5:k}
R_k\ge C\,\left(\frac{\|u_k \|_{\infty}\ }{\Big(\|f_k(\cdot,u_k) \|_{\frac{2N}{N+2}}\Big)^\frac{2N}{q(N+2)}\
M_k^{\ 1-\frac{2N}{q(N+2)}}}\right)^{1/\big(2-\frac{N}{q}\big)}.
\end{equation}

\bigskip

{\it Step 3. A lower bound for the $L^{2^*}$-norms.}

\medskip

\noindent Now, from definition of $B(x_k,R_k)$, 
\begin{equation*}
\int_{B(x_k,R_k)} |u_k|^{2^*} \ge  \left(\frac{1}{2}\|u_k\|_{\infty}\right)^{2^*}   \omega\,(R_k)^N,
\end{equation*}
where $\omega=\omega_N$ is the volume of the unit ball in $\R^N.$ 

\noindent Using   the inequality \eqref{Sobolev:5:k},  we deduce
\begin{align*}
\int_\Om |u_k|^{2^*}
&\ge  
C\,  \left( \frac{
(\|u_k\|_{\infty})^{1+2^*\left(\frac2{N}-\frac{1}{q}\right)}\ }
{\Big(\|f_k(\cdot,u_k) \|_{\frac{2N}{N+2}}\Big)^\frac{2N}{q(N+2)}\
M_k^{\ 1-\frac{2N}{q(N+2)}}}\right)^{\frac{1}{\frac2{N}-\frac{1}{q}}}.
\end{align*}
Denoting
$$
a=1+2^*\left(\frac2{N}-\frac{1}{q}\right),\qquad b=1-\frac{2N}{q(N+2)},
$$
observe that $\frac{a}{b}=2^*-1$. From (H2)$_k$,  
and due to $h_k$ is defined by \eqref{def:h} for $f=f_k$, we deduce
\begin{align*}
\|u_k\|_{2^*}^{2^*}
&\ge C\,  \Big(h_k\big(\|u_k\|_{\infty}\big)
\Big)^{\frac{\ 1-\frac{2N}{q(N+2)}}{\frac2{N}-\frac{1}{q}}}\
\frac1{\|f_k(\cdot,u_k) \|_{\frac{2N}{N+2}}}^{\frac{\frac{2N}{q(N+2)}}{\frac2{N}-\frac{1}{q}}}.
\end{align*}

Since the above, we can write  
\begin{equation*}\label{L:inf:h:4}
h_k\big(\|u_k\|_{\infty}\big)\le C\,\left[ \big(\|u_k\|_{2^*}\big)^{2^*\left(\frac2{N}-\frac{1}{q}\right)} \ \Big(\|f_k(\cdot,u_k) \|_{\frac{2N}{N+2}}\Big)^\frac{2N}{q(N+2)}\right]^\frac{1}{\ 1-\frac{2N}{q(N+2)}}.
\end{equation*}
Writing $\frac1{q}=\frac{1+\te}{N}$ with $\te\in(0,1)$ we obtain
\begin{equation*}\label{L:inf:h:6}
h_k\big(\|u_k\|_{\infty}\big)\le C\ \big(\|u_k\|_{2^*}\big)^{\frac{N+2}{N-2}\frac{2(1-\te)}{N-2\te}} \ \Big(\|f_k(\cdot,u_k) \|_{\frac{2N}{N+2}}\Big)^\frac{2(1+\te)}{N-2\te}.
\end{equation*}
Finally by elliptic regularity and Sobolev embedding,  
\begin{equation}\label{ell:reg:k:2}
\|u_k \|_{2^*}
\le C\ \|u_k \|_{H^{1}_0 (\Om)}
\le   C\
\|f_k(\cdot,u_k) \|_{\frac{2N}{N+2}},
\end{equation}
and we deduce that 
\begin{equation*}\label{L:inf:h:5}
h_k\big(\|u_k\|_{\infty}\big)\le C\, \Big(\|f_k(\cdot,u_k) \|_{\frac{2N}{N+2}}\Big)^{\frac4{N-2}},
\end{equation*}
ending  the proof. 
\end{proof}

\section{Approximation of  weak solutions. Sequences of BVP. Regularity theory of weak solutions.}
\label{sec:proof2}

In that Section, we  use families of BVP, and families of classical solutions  to approach some weak solutions, and prove Theorem \ref{th:incr:convex} and Theorem \ref{th:conv:seq}.

\medskip

The next Proposition  provides an approximation result for some  weak solutions to \eqref{eq:ell:pb} with subcritical nonlinearities. It states that there exist a family of BVP, and a family of solutions, uniformly bounded and convergent to a classical solution to \eqref{eq:ell:pb}. 

\begin{pro}\label{pro:conv:seq}
Let $f:\Omb \times[0,+\infty)\to [0,+\infty)$ be a  continuous function in both variables, satisfying {\rm (H1)} and \eqref{0:no:sol}. Assume that there exists $s_0>0$ such that  $f(x,\cdot)$ non-decreasing for all $s\ge s_0$, $x\in\Om$. 

Let $u^{*}\in H_0^1(\Om)$ be a   non-negative weak solution to \eqref{eq:ell:pb}.

\smallskip

Then, for some $\e_0>0$, there exist a family of non-linearities $\{f_\e \}_{\e \in (0,\e_0)}$ and a family of 
strong solutions  $\{u_\e\}_{\e\in (0,\e_0)}\subset W^{2,p}(\Om)\cap W_0^{1,p}(\Om)$ for any $p>1$, to
\begin{equation}\label{pde:e}
\left\{ 
\begin{array}{rcll}
- \De u_\e &=& f_\e(x,u_\e)  &\qquad\text{in}\ \Om,\\
u_\e &=& 0 &\qquad\text{on}\ \p \Om,
\end{array}
\right.
\end{equation}
such that $f_\e\le f$, $u_\e\le u^{*}$, and  $f_\e\uparrow f$ (pointwise in $\Omb\times\R$). Moreover, 
$u_\e\to \tilde{u}$ in $C^{1,\be}(\Omb))$ as $\e\to 0$, for any $\be<1$, and $\tilde{u}$ solves \eqref{eq:ell:pb}. 
\end{pro}

\begin{rem}
Once proved the above Proposition, there is still an open question: is $\tilde{u}= u^{*}$?. Theorem \ref{th:incr:convex} and Theorem \ref{th:conv:seq} shows two different ways to answer positively.
\end{rem}
\begin{proof}[Proof of Proposition \ref{pro:conv:seq}]
If $u^{*}\in L^\infty(\Om)$ the proof is easily achieved with $f_\e(x,s):=(1-\e)f(x,s)$.

Assume that $u^{*}\not\in L^\infty(\Om)$.

\bigskip

{\it Step 1. Construction of $f_\e\le f$.}

\medskip

\noindent Let us define
\begin{equation}\label{def:f:e}
f_\e(x,s):=\begin{cases}
(1-\e)f(x,s), & s\le 1/\e, \\
(1-\e)f(x,1/\e),\quad & s\ge 1/\e.
\end{cases}
\end{equation} 
Due to $f$ is non-negative 
and $f(x,\cdot)$ is non-decreasing for $s\ge s_0$, choosing $\e_0=1/s_0,$
$f_\e$ is a non decreasing family for any $\e\in(0,\e_0)$, 
$f_\e\le f$, and $f_\e\uparrow f$ (pointwise in $\Omb\times\R$) as $\e \downarrow0$.

Thanks to Beppo-Levi Theorem, 
and Lemma \ref{lem:w:H},
\begin{equation}\label{u:extrem}
f_\e(\cdot,u^{*}(\cdot))\to f(\cdot,u^{*}(\cdot))\qq{in} L^\frac{2N}{N+2}(\Om).
\end{equation}

Observe also that subcriticality (see \eqref{f:sub}) implies in particular the following
\begin{equation}\label{f:eps}
0\le f_\e(x,s)\le 
C_\e<+\infty,\qq{for all}s\ge 0,
\end{equation}
for each $\e\in(0,\e_0)$.

\bigskip

{\it Step 2. Construction  of $u_\e\le u^*$, strong solutions to \eqref{pde:e}$_\e$ in $W^{2,p}(\Om)\cap W_0^{1,p}(\Om)$ for any $p>1$.}

\medskip

\noindent Consider now the family of BVP's \eqref{pde:e}$_\e$. From \eqref{0:no:sol}, $0$ is a  subsolution to \eqref{pde:e}$_\e$ for any $\e\in(0,\e_0)$, and not a solution. 
On the other hand, since $f$ is non-negative and non-decreasing, $u^{*} $ is a 
weak supersolution to \eqref{pde:e}$_\e$ for any $\e\in(0,\e_0)$.

Consequently, there exist $L^1$-weak solutions $\underline{u}_\e \le  \overline{u}_\e$ of \eqref{pde:e}$_\e$ in $[0,u^{*}]$ such that any  $u_\e$ solution to \eqref{pde:e}$_\e$ in the interval $[0,u^{*}]$, satisfies
\begin{equation}\label{ineq:ue}
0\lneqq \underline{u}_\e \le u_\e\le \overline{u}_\e \le u^{*} \qq{a.e.}
\end{equation}
(see \cite[Theorem 1.1]{Montenegro-Ponce}), since \eqref{f:eps}, for each $\e>0,$ $f_\e$ is bounded,  hence $u_\e$ are strong solutions in $W^{2,p}(\Om)\cap W_0^{1,p}(\Om)$ for any $p>1$, and $u_\e\in C^{1,\be}(\Omb )$, see Proposition \ref{pro:classical}.

\bigskip

{\it Step 3. The family of solutions $u_\e$ is uniformly bounded in $C^{1,\be}(\Omb )$.}

\medskip

\noindent Fix $\e\in(0,\e_0)$. Since $f_\e\le f$, $f$ is non-negative, non-decreasing,  and  $u_\e\le u^*$
\begin{align*}\label{eq:ss:u:e}
\int_\Om |\nabla u_\e  |^2 \, dx
&=\int_\Om f_\e (x,u_\e)u_\e\le\int_\Om f (x,u_\e)u^{*}\\
&\le\int_\Om f (x,u^*)u^{*}
=\| u^{*}\|_{H_0^1(\Om)}^2\le C,
\end{align*}
where $C$ is only dependent on $f$ and $u^{*}$, and it is independent of $\e$. Now, Theorem \ref{th:apriori:cnys:k} implies that $\{u_\e\}$ are uniformly $L^\infty$ a priori bounded. By elliptic regularity  (see Proposition \ref{pro:classical}), there exists a uniform constant $C>0$ such that $\|u_\e\|_{C^{1,\be}(\Omb )}\le C$ for any $\be<1$. 

\bigskip

{\it Step 4. $u_\e\to \tilde{u}$  in $C^{1,\be}(\Omb))$ as $\e\to 0$, for any $\be<1$, and $\tilde{u}$ solves \eqref{eq:ell:pb}.}

\medskip

\noindent By compact embeddings and monotonicity, for any $\be'<\be<1$ the family $\{u_\e\}$, converges to $\tilde{u}$ in $C^{1,\be'}(\Omb )$ as $k\to\infty$, see  \cite[Lemma 6.36]{G-T}. 

Moreover, $\tilde{u}$ solves \eqref{eq:ell:pb}. Indeed, since $u_\e$ solves \eqref{pde:e}$_{\e}$, using the Lipschitzian property of $f$ on bounded intervals, and the uniform $L^\infty$ bound for $u_\e,$ and $\tilde{u}$, for any
$\varphi \in H_0^1(\Om)$:
\begin{align*}
&\left|\int_\Om \nabla \tilde{u} \nabla \varphi - f(x,\tilde{u})\varphi\, dx\right|
=\left|\int_\Om \nabla(u_\e-\tilde{u})  \nabla \varphi -  \big[f_{\e}(x,u_\e)-f(x,\tilde{u})\big]\varphi\Big) \, dx\right|\\
&\quad =\left|\int_\Om \nabla(u_\e-\tilde{u})  \nabla \varphi -  \big[f_{\e}(x,u_\e)-f(x,u_\e)+f(x,u_\e)-f(x,\tilde{u})\big]\varphi\Big) \, dx\right|
\\
&\quad \le C\Big(\|\nabla(u_\e-\tilde{u})\|_2 +  \e\big\|f(x,u_\e)\big\|_{\frac{2N}{N+2}}
+\big\|f(x,u_\e)-f(x,\tilde{u})\big\|_{\frac{2N}{N+2}}\Big)\|\vf\|_{H_0^1(\Om)}
\\
&\quad\to 0, \qq{as}k \to \infty, 
\end{align*}
ending the proof.
\end{proof}

\subsection{Proof of Theorem \ref{th:incr:convex}}
We extend the above result on smoothness  to the case of weak solutions to subcritical semilinear elliptic equations in any dimension. The question is now if for any positive semi-stable weak solution $u^*$, we can construct a sequence of BVP and a sequence of classical solutions convergent to  $u^*$.

\medskip

\begin{proof}[Proof of Theorem \ref{th:incr:convex}]
If $\sup_{\Om } u^*<+\infty$, then $u^{*}\in C^2(\Om)$ and the proof is finished. 

Assume  that 
\begin{equation}\label{u:inf}
\sup_{\Om } u^*=+\infty ,
\end{equation}
By Proposition \ref{pro:conv:seq}, there exists a family 
$u_\e\to \tilde{u}$ in $C^{2,\be}(\Omb))$ as $\e\to 0$, for any $\be<\al$, and $\tilde{u}$ solves \eqref{eq:ell:pb}. 
We will now prove that $\tilde{u}= u^*$. 

Assume by contradiction that $\tilde{u} \lneq u^*$. We observe that, by a density argument and Fatou's Lemma, the semi-stability inequality holds for $\vf\in H_0^1(\Om)$.
Testing it for $u^*$ with $u^* - \tilde{u}\gneq 0,$ we obtain
\begin{align*}
\int_\Om\big[f(x,u^*) - f(x,\tilde{u})\big](u^*- \tilde{u})\, dx
&= \int_\Om\big|\nabla (u^*- \tilde{u})\big| ^2 \, dx \\
&\ge  \int_\Om f_s(x,u^*)(u^*- \tilde{u}) ^2 \, dx.
\end{align*}
On the other hand,  by convexity
\begin{align*}
f(x,u^*) - f(x,\tilde{u})
&\le f_s(x,u^*)(u^*- \tilde{u}),
\end{align*}
this leads to
\begin{equation}\label{f:lin}
f(x,u^*) - f(x,\tilde{u})=f_s(x,u^*)(u^*- \tilde{u}) ,
\qq{a.e. in} \Om. 
\end{equation}

\bigskip

Due to  $\tilde{u}$ is a classical solution, there exists $\tilde{M}>0$ such that 
\begin{equation}\label{tilde:u:M}
0\le\tilde{u}\le \tilde{M}. 
\end{equation}
Hypothesis (H4) implies that given $\e_0=(c_0-1)/2>0$, there exists $s_0$ such that 
$$
\frac{sf_s(x,s)}{f(x,s)}
\ge \frac{c_0+1}{2}>1,\qq{for all} s\ge s_0,\ \text{a.e. in } \Om.
$$
Therefore, there exists $\de>0$ and $s_1>0$ such that 
$$
\frac{f_s(x,s)(s-\tilde{M})}{f(x,s)}
\ge \frac{(c_0+1)(s-\tilde{M})}{2s}\ge 1+\de,\qq{for all} s\ge s_1.
$$
Taking into account \eqref{u:inf}, there exists $\emptyset\not\equiv\om_1\subset \Om$ be such that $u^*(x)\ge s_1$ 
a.e. $x\in\om_1$.

From the above, and taking into account \eqref{f:lin}-\eqref{tilde:u:M},
we deduce
$$
1=\frac{f(x,\tilde{u}(x)) +f_s(x,u^*(x))(u^*(x)- \tilde{u}(x))}{f(x,u^*(x))}\ge 1+\de>1,\quad\text{a.e.}\ x\in\om_1 
$$
reaching a contradiction and concluding the proof.
\end{proof}

\subsection{Proof of Theorem \ref{th:conv:seq}}

The only difference with the proof of Theorem \ref{th:incr:convex} are the arguments involved in proving that $\tilde{u}= u^*$.

\medskip

\begin{proof}[Proof of Theorem \ref{th:conv:seq}]	
We apply Proposition \ref{pro:conv:seq} and get a sequence $u_k\to u$  in $C^{2,\be}(\Omb))$ as $k\to \infty$, for any $\be<\al$, and $u$ solves \eqref{eq:ell:pb}. 
We can repeat the argument for the family of maximal solutions  $\overline{u}_\e$ getting a subsequence  $\{\overline{u}_k\}$ convergent to $\overline{u}\le u^{*}$ in $C^{1,\be}(\Omb )$. Moreover   $\overline{u}\in C^{2,\be}(\Omb )$ solves \eqref{eq:ell:pb}. If $\overline{u}= u^{*}$, then the proof is finished.

Assume on the contrary that $\overline{u}\lneq u^{*}$. Fix $\e\in(0,\e_0)$ small enough. We now prove that  $\overline{u}_\e= u^{*}$, arguing on the contrary, and assuming that $\overline{u}_\e\lneq u^{*}$. Observe firstly that hypothesis (H5) imply  that there exists a positive constant $\de_0$ such that
\begin{equation}\label{fs:0}
f_s(x,s)\ge \la_1+\de_0,\qq{for any} x\in\Omb,\ s\ge 0,\\
\end{equation}
and consequently
\begin{equation}\label{f:0} 
f(x,s)\ge (\la_1+\de_0)\, s,\quad\text{for any}\ x\in\Omb,\ s\ge 0.
\end{equation}

Secondly, we see that $\overline{u}_\e+\de_1\phi_1$ is a subsolution to \eqref{pde:e}$_\e$ for any $\de_1>0$. From \eqref{fs:0}, and thanks to the mean value theorem, there exists a $\te=\te(x)$ such that
\begin{align*}
- \Delta(\overline{u}_\e+\de_1\phi_1\Big)
&=(1-\e)f(x,\overline{u}_\e)+\de_1\la_1\phi_1\le (1-\e)f(x,\overline{u}_\e+\de_1\phi)\\
&\iff \de_1\la_1\phi_1\le (1-\e)\big[f(x,\overline{u}_\e+\de_1\phi)-f(x,\overline{u}_\e)\big]
\\
&\iff \de_1\la_1\phi_1\le (1-\e)
f_s\big(x, \overline{u}_\e+\te\de_1\phi_1\big)\de_1\phi_1\\
&\iff \la_1\le (1-\e)
f_s\big(x,\overline{u}_\e+\te\de_1\phi_1\big) \ \checkmark .
\end{align*}

Thirdly, we check that  $\overline{u}_\e+\de_1\phi_1\le u^{*}$  for  $\de_1$ small enough.
From \eqref{f:0}, $f(x,\overline{u}_\e)\ge\la_1\overline{u}_\e$, and since  $\overline{u}_\e$ is a classical solution, there exists a $\de_2>0$ such that $\overline{u}_\e>\de_2\phi_1$.
Hence, for any $\varphi \in C^2(\Omb ),$ $\varphi\big|_{\p\Om}=0$, $\vf>0$, the following holds
\begin{align*}
\int_\Om \nabla\big(u^{*}&-\overline{u}_\e-\de_1\phi_1\big)\cdot\nabla \varphi =\int_\Om \Big[f(x,u^{*})-(1-\e)f(x,\overline{u}_\e)-\de_1\la_1\phi_1\Big]\vf\\
&= \int_\Om \Big[f_s\big(x,\te u^{*}+(1-\te)\overline{u}_\e\big)\, (u^{*}-\overline{u}_\e)
+\e f(x,\overline{u}_\e)-\de_1\la_1\phi_1\Big]\vf\\
&\ge \int_\Om \Big[\la_1\, (u^{*}-\overline{u}_\e) +\e\la_1 \overline{u}_\e-\de_1\la_1\phi_1\Big]\vf\\
&\ge \int_\Om \Big[\la_1\, (u^{*}-\overline{u}_\e) +\big(\e\de_2-\de_1\big)\la_1\phi_1\Big]\vf\ge 0,
\end{align*}
choosing $\de_1\le \e\de_2$. Therefore $\overline{u}_\e+\de_1\phi_1\le u^{*}$.
Consequently,  there exist a weak solution $\hat{u}_\e $ of \eqref{pde:e}$_\e$ in $[\overline{u}_\e+\de_1\phi_1,u^{*}]$, 
contradicting that $\overline{u}_\e$ is a maximal solution to \eqref{pde:e}$_\e$ in the interval $[0,u^{*}]$.

Consequently, $u_\e \uparrow u^{*}\in C^{2,\be}(\Omb )$, and finally, the elliptic regularity ends the proof.
\end{proof}

\bigskip

\section{Acknowledgments}

I would like to thank Professor Jos\'e Arrieta for helpful discussions.

The author was partially supported by Grant MTM2019-75465, MICINN, Spain and Grupo de Investigaci\'on CADEDIF 920894, UCM.

\def\cprime{$'$}

\end{document}